\documentclass[12pt,a4paper]{article}
\usepackage{amsmath}
\usepackage{amssymb}
\usepackage{enumerate}
\usepackage{t1enc}
\usepackage[latin1]{inputenc}
\usepackage[english]{babel}
\usepackage{enumerate}

\usepackage{mathtools}

\usepackage{amsfonts}
\usepackage{latexsym}
\usepackage{bm}
\newtheorem{theorem}{Theorem}[section]

\newtheorem{lemma}[theorem]{Lemma}
\newtheorem{cor}[theorem]{Corollary}

\begin{document}
\title{Isolation of $k$-cliques\medskip\medskip
}

\author{Peter Borg \\[2mm]
\normalsize Department of Mathematics \\
\normalsize Faculty of Science \\
\normalsize University of Malta\\
\normalsize Malta\\
\normalsize \texttt{peter.borg@um.edu.mt}
\and
Kurt Fenech \\ [2mm]
\normalsize Department of Mathematics \\
\normalsize Faculty of Science \\
\normalsize University of Malta\\
\normalsize Malta\\
\normalsize \texttt{kurt.fenech.10@um.edu.mt}\\ [5mm]
\and
Pawaton Kaemawichanurat 
\\ [2mm]
\normalsize Department of Mathematics \\
\normalsize King Mongkut's University of Technology Thonburi\\
\normalsize Thailand\\
\normalsize \texttt{pawaton.kae@kmutt.ac.th}
}

\date{}
\maketitle

\begin{abstract}
For any positive integer $k$ and any $n$-vertex graph $G$, let $\iota(G,k)$ denote the size of a smallest set $D$ of vertices of $G$ such that the graph obtained from $G$ by deleting the closed neighbourhood of $D$ contains no $k$-clique. Thus, $\iota(G,1)$ is the domination number of $G$. We prove that if $G$ is connected, then $\iota(G,k) \leq \frac{n}{k+1}$ unless $G$ is a $k$-clique or $k = 2$ and $G$ is a $5$-cycle. The bound is sharp. The case $k=1$ is a classical result of Ore, and the case $k=2$ is a recent result of Caro and Hansberg. Our result solves a problem of Caro and Hansberg. 
\end{abstract}

\section{Introduction}
Unless stated otherwise, we use small letters such as $x$ to denote non-negative integers or elements of a set, and capital letters such as $X$ to denote sets or graphs. The set of positive integers is denoted by $\mathbb{N}$. For $n \in \{0\} \cup \mathbb{N}$, the set $\{i \in \mathbb{N} \colon i \leq n\}$ is denoted by $[n]$. Note that $[0]$ is the empty set $\emptyset$. Arbitrary sets are assumed to be finite. For a set $X$, the set of $2$-element subsets of $X$ is denoted by ${X \choose 2}$ (that is, ${X \choose 2} = \{ \{x,y \} \colon x,y \in X, x \neq y \}$). 

If $Y$ is a subset of ${X \choose 2}$ and $G$ is the pair $(X,Y)$, then $G$ is called a \emph{graph}, $X$ is called the \emph{vertex set of $G$} and is denoted by $V(G)$, and $Y$ is called the \emph{edge set of $G$} and is denoted by $E(G)$. A \emph{vertex of $G$} is an element of $V(G)$, and an \emph{edge of $G$} is an element of $E(G)$. We call $G$ an \emph{$n$-vertex graph} if $|V(G)| = n$. We may represent an edge $\{v,w\}$ by $vw$. If $vw \in E(G)$, then 
we say that $w$ is a \emph{neighbour of $v$ in $G$} (and vice-versa).  
For $v \in V(G)$, $N_{G}(v)$ denotes the set of neighbours of $v$ in $G$, $N_{G}[v]$ denotes $N_{G}(v) \cup \{ v \}$, and $d_{G}(v)$ denotes $|N_{G} (v)|$ and is called the \emph{degree of $v$ in $G$}. 
For $S \subseteq V(G)$, 
$N_G[S]$ denotes $\bigcup_{v \in S} N_G[v]$ (the \emph{closed neighbourhood of $S$}), $G[S]$ denotes the graph $(S,E(G) \cap {S \choose 2})$ (the subgraph of $G$ \emph{induced by $S$}), and $G-S$ denotes $G[V(G) \backslash S]$ (the graph obtained from $G$ by \emph{deleting $S$}). We may abbreviate $G-\{v\}$ to $G-v$. Where no confusion arises, the subscript $G$ is omitted from any of the notation above that uses it; for example, $N_G(v)$ is abbreviated to $N(v)$.

If $G$ and $H$ are graphs, $f : V(H) \rightarrow V(G)$ is a bijection, and $E(G) = \{f(v)f(w) \colon vw \in E(H)\}$, then we say that $G$ is a \emph{copy of $H$}, and we write $G \simeq H$. Thus, a copy of $H$ is a graph obtained by relabeling the vertices of $H$.

For $n \geq 1$, the graphs $([n], {[n] \choose 2})$ and $([n], \{\{i,i+1\} \colon i \in [n-1]\})$ are denoted by $K_n$ and $P_n$, respectively. 
A copy of $K_n$ is called a \emph{complete graph} or an \emph{$n$-clique}. A copy of $P_n$ is called an \emph{$n$-path} or simply a \emph{path}. 

If $G$ and $H$ are graphs such that $V(H) \subseteq V(G)$ and $E(H) \subseteq E(G)$, then $H$ is called a \emph{subgraph of $G$}, and we say that \emph{$G$ contains $H$}.

If $\mathcal{F}$ is a set of graphs and $F$ is a copy of a graph in $\mathcal{F}$, then we call $F$ an \emph{$\mathcal{F}$-graph}. If $G$ is a graph and $D \subseteq V(G)$ such that $G-N[D]$ contains no $\mathcal{F}$-graph, then $D$ is called an \emph{$\mathcal{F}$-isolating set of $G$}. Let $\iota(G, \mathcal{F})$ denote the size of a smallest $\mathcal{F}$-isolating set of $G$. 
The study of isolating sets was introduced recently by Caro and Hansberg~\cite{CaHa17}. It is an appealing and natural generalization of the classical domination problem \cite{C, CH, HHS, HHS2, HL, HL2}. Indeed, $D$ is a $\{K_1\}$-isolating set of $G$ if and only if $D$ is a \emph{dominating set of $G$} (that is, $N[D] = V(G)$), so $\iota(G,\{K_1\})$ is the \emph{domination number of $G$} (the size of a smallest dominating set of $G$). In this paper, we obtain a sharp upper bound for $\iota(G, \{K_k\})$, and consequently we solve a problem of Caro and Hansberg \cite{CaHa17}.  

We call a subset $D$ of $V(G)$ a \emph{$k$-clique isolating set of $G$} if $G-N[D]$ contains no $k$-clique. We denote the size of a smallest $k$-clique isolating set of $G$ by $\iota(G,k)$. Thus, $\iota(G,k) = \iota(G, \{K_k\})$.

If $G_1, \dots, G_t$ are graphs such that $V(G_i) \cap V(G_j) = \emptyset$ for every $i,j \in [t]$ with $i \neq j$, then $G_1, \dots, G_t$ are \emph{vertex-disjoint}. A graph $G$ is \emph{connected} if, for every $v, w \in V(G)$, $G$ contains a path $P$ with $v, w \in V(P)$. A connected subgraph $H$ of $G$ is a \emph{component of $G$} if, for each connected subgraph $K$ of $G$ with $K \neq H$, $H$ is not a subgraph of $K$. Clearly, two distinct components of $G$ are vertex-disjoint.  

For $n, k \in \mathbb{N}$, let $a_{n,k} = \left \lfloor \frac{n}{k+1} \right \rfloor$ and $b_{n,k} = n - ka_{n,k}$. We have $a_{n,k} \leq b_{n,k} \leq a_{n,k} + k$. If $n \leq k$, then let $B_{n,k} = P_n$. If $n \geq k+1$, then let $F_1, \dots, F_{a_{n,k}}$ be copies of $K_k$ such that $P_{b_{n,k}}, F_1, \dots, F_{a_{n,k}}$ are vertex-disjoint, and let $B_{n,k}$ be the connected $n$-vertex graph given by 
\[B_{n,k} = \left( V(P_{b_{n,k}}) \cup \bigcup_{i=1}^{a_{n,k}} V(F_i), \, E(P_{b_{n,k}}) \cup \{iv \colon i \in [a_{n,k}], v \in V(F_i)\} \cup \bigcup_{i=1}^{a_{n,k}} E(F_i) \right).\] 
Thus, $B_{n,k}$ is the graph obtained by taking $P_{b_{n,k}}, F_1, \dots, F_{a_{n,k}}$ and joining $i$ (a vertex of $P_{b_{n,k}}$) to each vertex of $F_i$ for each $i \in [a_{n,k}]$.

For $n, k \in \mathbb{N}$ with $k \neq 2$, let 
\[\iota(n,k) = \max \{\iota(G,k) \colon G \mbox{ is a connected graph}, V(G) = [n], G \not\simeq K_k\}.\] 
For $n \in \mathbb{N}$, let 
\[\iota(n,2) = \max \{\iota(G,2) \colon G \mbox{ is a connected graph}, V(G) = [n], G \not\simeq K_2, G \not\simeq C_5\}.\]
In Section~\ref{Proofsection}, we prove the following result.

\begin{theorem} \label{mainresult1}
If $G$ is a connected $n$-vertex graph, then, unless $G$ is a $k$-clique or $k = 2$ and $G$ is a $5$-cycle,
\[ \iota(G,k) \leq \frac{n}{k+1}. \]
Consequently, for any $k \geq 1$ and $n \geq 3$, 
\[\iota(n,k) = \iota(B_{n,k}, k) = \left \lfloor \frac{n}{k+1} \right \rfloor .\]
\end{theorem}

A classical result of Ore \cite{Ore} is that the domination number of a graph $G$ with $\min\{d(v) \colon v \in V(G)\} \geq 1$ is at most $\frac{n}{2}$ (see \cite{HHS}). Since the domination number is $\iota(G,1)$, it follows by Lemma~\ref{lemmacomp} in Section~\ref{Proofsection} that Ore's result is equivalent to the bound in Theorem~\ref{mainresult1} for $k = 1$. The case $k = 2$ is also particularly interesting; while deleting the closed neighbourhood of a $\{K_1\}$-isolating set yields the graph with no vertices, deleting the closed neighbourhood of a $\{K_2\}$-isolating set yields a graph with no edges. In \cite{CaHa17}, Caro and Hansberg proved Theorem~\ref{mainresult1} for $k = 2$, using a different argument. Consequently, they established that $\frac{1}{k+1} \leq \limsup_{n \rightarrow \infty} \frac{\iota(n,k)}{n} \leq \frac{1}{3}$. In 
the same paper, they asked for the value of $\limsup_{n \rightarrow \infty} \frac{\iota(n,k)}{n}$. The answer is given by Theorem~\ref{mainresult1}.

\begin{cor} For any $k \geq 1$,
\[\limsup_{n \rightarrow \infty} \frac{\iota(n,k)}{n} = \frac{1}{k+1}.\]
\end{cor}
\textbf{Proof.} By Theorem~\ref{mainresult1}, for any $n \geq 3$, $\frac{1}{k+1} - \frac{k}{(k+1)n} = \frac{1}{n}\left( \frac{n-k}{k+1} \right) \leq \frac{\iota(n,k)}{n} \leq \frac{1}{k+1}$, and, if $n$ is a multiple of $k+1$, then $\frac{\iota(n,k)}{n} = \frac{1}{k+1}$. Thus, $\lim_{n \rightarrow \infty} \sup \left\{\frac{\iota(p,k)}{p} \colon p \geq n \right\}  = \lim_{n \rightarrow \infty} \frac{1}{k+1} = \frac{1}{k+1}$.~\hfill{$\Box$}

\section{Proof of Theorem~\ref{mainresult1}} \label{Proofsection}

In this section, we prove Theorem~\ref{mainresult1}. We start with two lemmas that will be used repeatedly. 

If a graph $G$ contains a $k$-clique $H$, then we call $H$ a \emph{$k$-clique of $G$}. We denote the set $\{V(H) \colon H \mbox{ is a $k$-clique of } G\}$ by $\mathcal{C}_k(G)$. 

\begin{lemma} \label{lemma}
If $v$ is a vertex of a graph $G$, then $\iota(G,k) \leq 1 + \iota(G-N_G[v],k)$. 
\end{lemma}
\textbf{Proof.} Let $D$ be a $k$-clique isolating set of $G-N_G[v]$ of size $\iota(G-N_G[v],k)$. Clearly, $C \cap N_G[v] \neq \emptyset$ for each $C \in \mathcal{C}_k(G) \backslash \mathcal{C}_k(G-N_G[v])$. Thus, $D \cup \{v \}$ is a $k$-clique isolating set of $G$. The result follows. ~\hfill{$\Box$}

\begin{lemma} \label{lemmacomp}
If $G_1, \dots, G_r$ are the distinct components of a graph $G$, then $\iota(G,k) = \sum_{i=1}^r \iota(G_i,k)$. 
\end{lemma}
\textbf{Proof.} For each $i \in [r]$, let $D_i$ be a smallest $k$-clique isolating set of $G_i$. Then, $\bigcup_{i=1}^r D_i$ is a $k$-clique isolating set of $G$. Thus, $\iota(G,k) \leq \sum_{i=1}^r |D_i| = \sum_{i=1}^r \iota(G_i,k)$. Let $D$ be a smallest $k$-clique isolating set of $G$. For each $i \in [r]$, $D \cap V(G_i)$ is a $k$-clique isolating set of $G_i$, so $|D_i| \leq |D \cap V(G_i)|$. We have $\sum_{i=1}^r \iota(G_i,k) = \sum_{i=1}^r |D_i| \leq \sum_{i=1}^r |D \cap V(G_i)| = |D| = \iota(G,k)$. The result follows.~\hfill{$\Box$}
\\
\\
\textbf{Proof of Theorem~\ref{mainresult1}.} We use induction on $n$. If $G$ is a $k$-clique, then $\iota(G) = 1 = \frac{n+1}{k+1}$. If $k = 2$ and $G$ is a $5$-cycle, then $\iota(G) = 2 = \frac{n+1}{k+1}$. Suppose that $G$ is not a $k$-clique and that, if $k = 2$, then $G$ is not a $5$-cycle. Suppose $n \leq 2$. If $k \geq 3$, then $\iota(G) = 0$. If $k = 2$, then $G \simeq K_1$, so $\iota(G) = 0$. If $k = 1$, then $G \simeq K_2$, so $\iota(G) = 1 = \frac{n}{k+1}$. Now suppose $n \geq 3$. If $\mathcal{C}_k(G) = \emptyset$, then $\iota_k(G) = 0$. Suppose $\mathcal{C}_k(G) \neq \emptyset$. 
Let $C \in \mathcal{C}_k(G)$. Since $G$ is connected and $G$ is not a $k$-clique, there exists some $v \in C$ such that $N[v] \backslash C \neq \emptyset$. Thus, $|N[v]| \geq k+1$ as $C \subset N[v]$. If $V(G) = N[v]$, then $\{v\}$ is a $k$-clique isolating set of $G$, so $\iota(G) = 1 \leq \frac{n}{k+1}$. Suppose $V(G) \neq N[v]$. Let $G' = G-N[v]$ and $n' = |V(G')|$. Then, 
\[n \geq n' + k + 1\] 
and $V(G') \neq \emptyset$. Let $\mathcal{H}$ be the set of components of $G'$. If $k \neq 2$, then let $\mathcal{H}' = \{H \in \mathcal{H} \colon H \simeq K_k\}$. If $k = 2$, then let $\mathcal{H}' = \{H \in \mathcal{H} \colon H \simeq K_k\mbox{ or } H \simeq C_5\}$. By the induction hypothesis, $\iota(H,k) \leq \frac{|V(H)|}{k+1}$ for each $H \in \mathcal{H} \backslash \mathcal{H}'$. If $\mathcal{H}' = \emptyset$, then, by Lemmas~\ref{lemma} and \ref{lemmacomp},
\begin{equation}
\iota(G,k) \leq 1 + \iota(G',k) = 1 + \sum_{H \in \mathcal{H}} \iota(H,k) \leq 1 + \sum_{H \in \mathcal{H}} \frac{|V(H)|}{k+1} = \frac{k+1 + n'}{k+1} \leq \frac{n}{k+1}. \nonumber
\end{equation}

Suppose $\mathcal{H}' \neq \emptyset$. For any $H \in \mathcal{H}$ and any $x \in N(v)$, we say that $H$ is \emph{linked to $x$} if $xy \in E(G)$ for some $y \in V(H)$. Since $G$ is connected, each member of $\mathcal{H}$ is linked to at least one member of $N(v)$. One of Case 1 and Case 2 below holds.
\\
\\
\emph{Case 1: For each $H \in \mathcal{H}'$, $H$ is linked to at least two members of $N(v)$.} Let $H' \in \mathcal{H}'$ and $x \in N(v)$ such that $H'$ is linked to $x$. Let $\mathcal{H}_x$ be the set of members of $\mathcal{H}$ that are linked to $x$ only. Then, 
\[\mathcal{H}_x \subseteq \mathcal{H} \backslash \mathcal{H}',\] 
and hence, by the induction hypothesis, each member $H$ of $\mathcal{H}_x$ has a $k$-clique isolating set $D_H$ with $|D_H| \leq \frac{|V(H)|}{k+1}$. 

Let $X = \{x\} \cup V(H')$ and $G^* = G - X$. Then, $G^*$ has a component $G^*_v$ with $N[v] \backslash \{x\} \subseteq V(G^*_v)$, and the other components of $G^*$ are the members of $\mathcal{H}_x$. Let $D^*_v$ be a $k$-clique isolating set of $G^*_v$ of size $\iota(G^*_v,k)$. Since $H'$ is linked to $x$, $xy \in E(G)$ for some $y \in V(H')$. If $H'$ is a $k$-clique, then let $D' = \{y\}$. If $k = 2$ and $H'$ is a $5$-cycle, then let $y'$ be one of the two vertices in $V(H') \backslash N_{H'}[y]$, and let $D' = \{y, y'\}$. We have $X \subseteq N[D']$ and $|D'| = \frac{|X|}{k+1}$. Let $D = D' \cup D^*_v \cup \bigcup_{H \in \mathcal{H}_x} D_H$. Since the components of $G^*$ are $G^*_v$ and the members of $\mathcal{H}_x$, we have $V(G) = X \cup V(G^*_v) \cup \bigcup_{H \in \mathcal{H}_x} V(H)$, and, since $X \subseteq N[D']$, $D$ is a $k$-clique isolating set of $G$. Thus,
\begin{equation} \iota(G,k) \leq |D| = |D^*_v| + |D'| + \sum_{H \in \mathcal{H}_x} |D_H| \leq |D^*_v| + \frac{|X|}{k+1} + \sum_{H \in \mathcal{H}_x} \frac{|V(H)|}{k+1}. \label{main.1}
\end{equation}
%
\emph{Subcase 1.1: $G^*_v$ is neither a $k$-clique nor a $5$-cycle.} Then, $|D^*_v| \leq \frac{|V(G^*_v)|}{k+1}$ by the induction hypothesis. By (\ref{main.1}), $\iota(G,k) \leq \frac{1}{k+1} \left( |V(G^*_v)| + |X| + \sum_{H \in \mathcal{H}_x} |V(H)| \right) = \frac{n}{k+1}$.\medskip
\\
\emph{Subcase 1.2: $G^*_v$ is a $k$-clique.} Since $|N[v]| \geq k+1$ and $N[v] \backslash \{x\} \subseteq V(G^*_v)$, we have $V(G^*_v) = N[v] \backslash \{x\}$. If $H'$ is a $k$-clique, then let $X' = \{y\}$ and $D'' = \{x\}$. If $k = 2$ and $H'$ is a $5$-cycle, then let $X'$ be the set whose members are $y$, $y'$, and the two neighbours of $y'$ in $H'$, and let $D'' = \{x, y'\}$. Let $Y = (X \cup V(G^*_v)) \backslash (\{v,x\} \cup X')$. Let $G_Y = G - (\{v,x\} \cup X')$. Then, the components of $G_Y$ are the components of $G[Y]$ and the members of $\mathcal{H}_x$. 

If $G[Y]$ has no $k$-clique, then, since $\{v,x\} \cup X' \subseteq N[D'']$, $D'' \cup \bigcup_{H \in \mathcal{H}_x} D_H$ is a $k$-clique isolating set of $G$, and hence 
\[\iota(G,k) \leq |D''| + \sum_{H \in \mathcal{H}_x} |D_H| < \frac{|X \cup V(G^*_v)|}{k+1} + \sum_{H \in \mathcal{H}_x} \frac{|V(H)|}{k+1} = \frac{n}{k+1}.\]
This is the case if $k = 1$ as we then have $Y = \emptyset$. 

Suppose that $k \geq 2$ and $G[Y]$ has a $k$-clique $C_Y$. We have
\begin{equation} V(C_Y) \subseteq (V(G^*_v) \backslash \{v\}) \cup (V(H') \backslash X'). \label{main.2}
\end{equation}
Thus, $|V(C_Y) \cap V(G^*_v)| = |V(C_Y) \backslash (V(H') \backslash X')| \geq k - (k-1) = 1$ and $|V(C_Y) \cap V(H')| = |V(C_Y) \backslash (V(G^*_v) \backslash \{v\})| \geq k - (k-1) = 1$. Let $z \in V(C_Y) \cap V(G^*_v)$ and $Z = V(G^*_v) \cup V(C_Y)$. Since $z$ is a vertex of each of the $k$-cliques $G^*_v$ and $C_Y$, 
\begin{equation} Z \subseteq N[z]. \label{main.3}
\end{equation}
We have 
\begin{equation} |Z| = |V(G^*_v)| + |V(C_Y) \backslash V(G^*_v)| = k + |V(C_Y) \cap V(H')| \geq k+1. \label{main.4}
\end{equation}
Let $G_Z = G - Z$. Then, $V(G_Z) = \{x\} \cup (V(H') \backslash V(C_Y)) \cup \bigcup_{H \in \mathcal{H}_x} V(H)$. We have that the components of $G_Z - x$ are $G_Z[V(H') \backslash V(C_Y)]$ (which is a clique or a path, depending on whether $H'$ a $k$-clique or a $5$-cycle) and the members of $\mathcal{H}_x$, $y \in V(H') \backslash V(C_Y)$ (by (\ref{main.2})), $y \in N_{G_Z}[x]$, and, by the definition of $\mathcal{H}_x$, $N_{G_Z}(x) \cap V(H) \neq \emptyset$ for each $H \in \mathcal{H}_x$. Thus, $G_Z$ is connected, and, if $\mathcal{H}_x \neq \emptyset$, then $G_Z$ is neither a clique nor a $5$-cycle. 

Suppose $\mathcal{H}_x \neq \emptyset$. By the induction hypothesis, $\iota(G_Z,k) \leq \frac{|V(G_Z)|}{k+1}$. Let $D_{G_Z}$ be a $k$-clique isolating set of $G_Z$ of size $\iota(G_Z,k)$. By (\ref{main.3}), $\{z\} \cup D_{G_Z}$ is a $k$-clique isolating set of $G$. Thus, $\iota(G,k) \leq 1 + \iota_k(G_Z) \leq 1 + \frac{|V(G_Z)|}{k+1}$, and hence, by (\ref{main.4}), $\iota(G,k) \leq \frac{|Z|}{k+1} + \frac{|V(G_Z)|}{k+1} = \frac{n}{k+1}$.

Now suppose $\mathcal{H}_x = \emptyset$. Then, $G^* = G^*_v$, so $V(G) = V(G^*_v) \cup \{x\} \cup V(H')$. Recall that either $H'$ is a $k$-clique or $k=2$ and $H'$ is a $5$-cycle. 

Suppose that $H'$ is a $k$-clique. Then, $n = 2k+1$. By (\ref{main.3}), $|V(G-N[z])| \leq |V(G-Z)| = n - |Z| = 2k+1 - |Z|$. Suppose $|Z| \geq k+2$. Then, $|V(G-N[z])| \leq k-1$, and hence $\{z\}$ is a $k$-clique isolating set of $G$. Thus, $\iota(G,k) = 1 < \frac{n}{k+1}$. Now suppose $|Z| \leq k+1$. Then, by (\ref{main.4}), $|Z| = k+1$ and $|V(C_Y) \cap V(H')| = 1$. 
Let $z'$ be the element of $V(C_Y) \cap V(H')$, and let $Z' = V(C_Y) \cup V(H')$. Since $z'$ is a vertex of each of the $k$-cliques $C_Y$ and $H'$, $Z' \subseteq N[z']$. We have $|Z'| = |V(C_Y)| + |V(H')| - |V(C_Y) \cap V(H')| = 2k-1$ and $|V(G-N[z'])| \leq |V(G-Z')| = n - |Z'| = (2k+1) - (2k-1) = 2$. If $k \geq 3$, then $\{z'\}$ is a $k$-clique isolating set of $G$, and hence $\iota(G) = 1 < \frac{n}{k+1}$. Suppose $k = 2$. Then, $H'$, $G^*_v$, and $C_Y$ are the $2$-cliques with vertex sets $\{y,z'\}$, $\{v,z\}$, and $\{z,z'\}$, respectively. Thus, $V(G) = \{v, z, z', y, x\}$, and $G$ contains the $5$-cycle with edge set $\{xv, vz, zz', z'y, yx\}$. Since $G$ is not a $5$-cycle, $d(w) \geq 3$ for some $w \in V(G)$. Since $|V(G - N[w])| = 5 - |N[w]| \leq 1$, $\{w\}$ is a $k$-clique isolating set of $G$, and hence $\iota(G,k) = 1 < \frac{5}{3} = \frac{n}{k+1}$.

Now suppose that $k=2$ and $H'$ is a $5$-cycle. Then, $V(G^*_v) = \{v, z\}$ and $E(H') = \{yy_1, y_1y_2, y_2y_3, y_3y_4, y_4y\}$ for some $y_1, y_2, y_3, y_4 \in V(G)$. Recall that $|V(C_Y) \cap V(H')| \geq 1$. Let $z' \in V(C_Y) \cap V(H')$. Since $z$ and $z'$ are vertices of $C_Y$, $zz' \in E(G)$. We have $V(G) = \{v, z, x, y, y_1, y_2, y_3, y_4\}$, $N(v) = \{x,z\}$, $z' \in \{y_1, y_2, y_3, y_4\}$ (as $y \notin V(C_Y)$ by (\ref{main.2})), and $\{vx, vz, xy, zz'\} \cup E(H') \subseteq E(G)$. If $z'$ is $y_1$ or $y_4$, then $V(G - N[\{y,z'\}])$ is $\{v, y_3\}$ or $\{v, y_2\}$. If $z'$ is $y_2$ or $y_3$, then $V(G - N[\{y,z'\}]) = \{v\}$. Thus, $\{y,z'\}$ is a $k$-clique reducing set of $G$, and hence $\iota(G,k) = 2 < \frac{8}{3} = \frac{n}{k+1}$.\medskip
\\
\emph{Subcase 1.3: $G^*_v$ is a $5$-cycle.} If $k \neq 2$, then the result follows as in Subcase~1.1. Suppose $k = 2$. We have $E(G^*_v) = \{vv_1, v_1v_2, v_2v_3, v_3v_4, v_4v\}$ for some $v_1, v_2, v_3, v_4 \in V(G)$. Let $Y = \{v_2, v_3, v_4\}$. Recall that the components of $G^*$ are $G^*_v$ and the members of $\mathcal{H}_x$. Thus, $G-Y$ is connected and $V(G-Y) = \{v, v_1, x\} \cup V(H') \cup \bigcup_{H \in \mathcal{H}_x} V(H)$. 

Suppose that $G-Y$ is not a $5$-cycle. By the induction hypothesis, $G-Y$ has a $k$-clique isolating set $D$ with $|D| \leq \frac{|V(G-Y)|}{k+1} = \frac{n-3}{3} = \frac{n}{3}-1$. Since $Y \subseteq N[v_3]$, $\{v_3\} \cup D$ is a $k$-clique isolating set of $G$, so $\iota(G,k) \leq \frac{n}{3} = \frac{n}{k+1}$.

Now suppose that $G-Y$ is a $5$-cycle. Then, $H'$ is a $2$-clique and $V(G-Y) = \{v,v_1,x,y,z\}$, where $\{z\} = V(H') \backslash \{y\}$. Since $v_1v, vx, xy, yz \in E(G-Y)$ and $G-Y$ is a $5$-cycle, $E(G-Y) = \{v_1v, vx, xy, yz, zv_1\}$. We have $V(G - N[\{v, v_1\}]) \subseteq \{v_3, y\}$. If $v_3y \notin E(G)$, then $\{v,v_1\}$ is a $k$-clique isolating set of $G$. If $v_3y \in E(G)$, then $V(G - N[v, v_3]) \subseteq \{z\}$, 
so $\{v,v_3\}$ is a $k$-clique isolating set of $G$. Therefore, $\iota(G,k) = 2 < \frac{8}{3} = \frac{n}{k+1}$.
\\
\\
\emph{Case 2: For some $x \in N(v)$ and some $H' \in \mathcal{H}'$, $H'$ is linked to $x$ only.} Let $\mathcal{H}_1 = \{H \in \mathcal{H}' \colon H$ is linked to $x$ only$\}$ and $\mathcal{H}_2 = \{H \in \mathcal{H} \backslash \mathcal{H}' \colon H$ is linked to $x$ only$\}$. Let $h_1 = |\mathcal{H}_1|$ and $h_2 = |\mathcal{H}_2|$. Since $H' \in \mathcal{H}_1$, $h_1 \geq 1$. For each $H \in \mathcal{H}_1$, $y_H \in N(x)$ for some $y_H \in V(H)$. Let $X = \{x\} \cup \bigcup_{H \in \mathcal{H}_1} V(H)$.  

For each $k$-clique $H \in \mathcal{H}_1$, let $D_H = \{x\}$. If $k = 2$, then, for each $5$-cycle $H \in \mathcal{H}_1$, let $y_H'$ be one of the two vertices in $V(H) \backslash N_{H}[y_H]$, and let $D_H = \{x, y_H'\}$. Let $D_X = \bigcup_{H \in \mathcal{H}_1} D_H$. Then, $D_X$ is a $k$-clique isolating set of $G[X]$. If $k \neq 2$, then $D_X = \{x\}$, so $|D_X| = 1 \leq \frac{1 + k|\mathcal{H}_1|}{k+1} = \frac{|X|}{k+1}$. If $k = 2$ and we let $h_1' = |\{H \in \mathcal{H}_1 \colon H \simeq C_5\}|$, then $|D_X| = 1 + h_1' \leq \frac{1 + 5h_1' + 2(h_1-h_1')}{3} = \frac{|X|}{k+1}$. 

Let $G^* = G-X$. Then, $G^*$ has a component $G^*_v$ with $N[v] \backslash \{x\} \subseteq V(G^*_v)$, and the other components of $G^*$ are the members of $\mathcal{H}_2$. By the induction hypothesis, $\iota(H,k) \leq \frac{|V(H)|}{k+1}$ for each $H \in \mathcal{H}_2$. For each $H \in \mathcal{H}_2$, let $D_H$ be a $k$-clique isolating set of $H$ of size $\iota(H,k)$. 

If $G^*_v$ is a $k$-clique, then let $D^*_v = \{x\}$. If $k=2$ and $G^*_v$ is a $5$-cycle, then let $v'$ be one of the two vertices in $V(G^*_v) \backslash N_{G^*_v}[v]$, and let $D^*_v = \{x, v'\}$. If neither $G^*_v$ is a $k$-clique nor $k=2$ and $G^*_v$ is a $5$-cycle, then, by the induction hypothesis, $G^*_v$ has a $k$-clique isolating set $D^*_v$ with $|D^*_v| \leq \frac{|V(G^*_v)|}{k+1}$. 

Let $D = D^*_v \cup D_X \cup \bigcup_{H \in \mathcal{H}_2} D_H$. By the definition of $\mathcal{H}_1$ and $\mathcal{H}_2$, the components of $G-x$ are $G^*_v$ and the members of $\mathcal{H}_1 \cup \mathcal{H}_2$. Thus, $D$ is a $k$-clique isolating set of $G$ since $x \in D$, $v \in V(G^*_v) \cap N[x]$, and $D_X$ is a $k$-clique isolating set of $G[X]$. Let $D' = D_X \cup \bigcup_{H \in \mathcal{H}_2} D_H$ and $n^* = |V(G^*_v)|$. We have 
\[|D'| = |D_X| + \sum_{H \in \mathcal{H}_2} |D_H| \leq \frac{|X|}{k+1} + \sum_{H \in \mathcal{H}_2} \frac{|V(H)|}{k+1} = \frac{n-n^*}{k+1}.\] 
If $G^*_v$ is a $k$-clique, then $|D| = |D'| < \frac{n}{k+1}$. If $k=2$ and $G^*_v$ is a $5$-cycle, then 
\[|D| = 1 + |D'| \leq 1 + \frac{n-n^*}{k+1} = 1 + \frac{n-5}{3} < \frac{n}{k+1}.\]
If neither $G^*_v$ is a $k$-clique nor $k=2$ and $G^*_v$ is a $5$-cycle, then $|D| = |D^*_v| + |D'| \leq \frac{n^*}{k+1} + \frac{n-n^*}{k+1} = \frac{n}{k+1}$.~\hfill{$\Box$}

\end{document}